\documentclass[11pt]{article} \topmargin=-0.5cm
\usepackage{amssymb}

   \csname @addtoreset\endcsname{equation}{section}
    \csname @addtoreset\endcsname{figure}{section}
    \csname @addtoreset\endcsname{table}{section}

\newcounter{thanksnum}
\def\thanksnumber#1
{\setcounter{thanksnum}{\value{footnote}}\setcounter{footnote}{#1}%
                     \addtocounter{footnote}{-1}\footnotemark
                     \setcounter{footnote}{\value{thanksnum}}}
\def\newtheoremz#1{\@ifnextchar[{\@othmz{#1}}{\@nthmz{#1}}}

\def\@nthmz#1#2{%
\@ifnextchar[{\@xnthmz{#1}{#2}}{\@ynthmz{#1}{#2}}}

\def\@xnthmz#1#2[#3]{\expandafter\@ifdefinable\csname #1\endcsname
{\@definecounter{#1}\@addtoreset{#1}{#3}%
\expandafter\xdef\csname the#1\endcsname{\expandafter\noexpand
  \csname the#3\endcsname \@thmcountersepz \@thmcounterz{#1}}%
\global\@namedef{#1}{\@thmz{#1}{#2}}\global\@namedef{end#1}{\@endtheoremz}}}

\def\@ynthmz#1#2{\expandafter\@ifdefinable\csname #1\endcsname
{\@definecounter{#1}%
\expandafter\xdef\csname the#1\endcsname{\@thmcounterz{#1}}%
\global\@namedef{#1}{\@thm{#1}{#2}}\global\@namedef{end#1}{\@endtheoremz}}}

\def\@othmz#1[#2]#3{\expandafter\@ifdefinable\csname #1\endcsname
  {\global\@namedef{the#1}{\@nameuse{the#2}}%
\global\@namedef{#1}{\@thmz{#2}{#3}}%
\global\@namedef{end#1}{\@endtheoremz}}}

\def\@thmz#1#2{\refstepcounter
    {#1}\@ifnextchar[{\@ythmz{#1}{#2}}{\@xthmz{#1}{#2}}}

\def\@xthmz#1#2{\@begintheoremz{#2}{\csname the#1\endcsname}\ignorespaces}
\def\@ythmz#1#2[#3]{\@opargbegintheoremz{#2}{\csname
       the#1\endcsname}{#3}\ignorespaces}

\def\@thmcounterz#1{\noexpand\arabic{#1}}
\def\@thmcountersepz{.}
\def\@begintheoremz#1#2{ \trivlist \item[\hskip \labelsep{\bf #1\ #2}]}
\def\@opargbegintheoremz#1#2#3{ \trivlist
      \item[\hskip \labelsep{\bf #1\ #2\ (#3)}]}
\def\@endtheoremz{\endtrivlist}

\newtheorem{theorem}{Theorem}[section]
\newtheorem{lemma}{Lemma}[section]

\newtheorem{proposition}{Proposition}[section]

\newtheorem{condition}{Condition}[section]
\newtheorem{definition}{Definition}[section]

\newtheorem{example}{Example}[section]

\def\e{\varepsilon}

\def\defi{\stackrel{{\scriptscriptstyle \Delta}}{=}}

\def\a{\alpha}
\def\d{\delta}
\def\o{\omega}
\def\O{\Omega}

\def\F{{\cal F}}
\def\w{\widehat}
\def\Ind{{\mathbb{I}}}

\def\esssup{\mathop{\rm ess\, sup}}

\def\R{{\bf R}}
\def\E{{\bf E}}
\def\P{{\bf P}}

\def\L{L}

\def\b{\beta}
\def\s{\delta}
\def\g{\gamma}

\def\W{{\cal W}^*}
\def\ww{\widetilde}
\def\X{{\cal X}}
\def\t{\theta}
\def\oo{\bar}
\def\s{\sigma}
\def\p{\partial}
\def\G{\Gamma}

\def\V{{\cal V}}
\def\A{{\cal A}}

\def\L{{\cal L}}

\def\TT{{\cal T}}
\newcommand{\be}{\begin{equation}}
\newcommand{\ee}{\end{equation}}
\newcommand{\bd}{\begin{displaymath}}
\newcommand{\ed}{\end{displaymath}}
\newcommand{\ba}{\begin{array}{ll}}
\newcommand{\ea}{\end{array}}
\newcommand{\baa}{\begin{eqnarray}}
\newcommand{\eaa}{\end{eqnarray}}
\newcommand{\baaa}{\begin{eqnarray*}}
\newcommand{\eaaa}{\end{eqnarray*}}   \font\sm=cmr10
\def\ww{\tilde}
\def\W{{\cal W}}

\date{Submitted: November 7, 2012. Revised: July 31, 2013}

\def\Q{{\cal Q}}
\def\CC{{\cal C}}
\def\iindex{}
\title{Backward SPDEs with non-local  in time and space boundary
conditions}% Arxiv: 1211.1460 
\author{
Nikolai Dokuchaev\\
 {\sm Department of Mathematics \& Statistics, Curtin
University,}\\ {\sm  GPO Box U1987, Perth, 6845 Western Australia} }
\begin{document}
\maketitle
\begin{abstract}
We study linear  backward stochastic partial differential equations
of parabolic type with special boundary condition that connect the
terminal value of the solution with a functional over the entire past
solution. Uniqueness, solvability and regularity results for the
solutions are obtained.
\end{abstract}
%{\it Abbreviated tittle: }
\section{Introduction}
Partial differential equations and  stochastic partial
differential equations (SPDEs) have fundamental significance for
natural sciences, and various boundary value problems for them
were widely studied.
 Usually,  well-posedness of a  boundary value depends on
the choice of the boundary value conditions.
\par
 Boundary value problems for SPDEs are well studied  in the existing literature
 for the case
 of  forward parabolic Ito equations with the  Cauchy condition at
initial time (see, e.g., Al\'os et al (1999), Bally {\it et al}
(1994),  Da Prato and Tubaro (1996), Gy\"ongy (1998), Krylov (1999),
Maslowski (1995), Pardoux (1993),
 Rozovskii (1990), Walsh (1986), Zhou (1992),
and the bibliography there).  Many results have been also obtained
for the backward
 parabolic Ito equations with  Cauchy
 condition   at terminal time, as well as for  pairs of   forward and
backward
 equations with separate  Cauchy conditions at initial time and
 the terminal time respectively; see, e.g., Yong and Zhou
 (1999), and the author's papers (1992), (2005),(2011), (2012a).  Note that a
 backward SPDE cannot be transformed into a forward
  equation by a simple time
 change, unlike as for the case of deterministic equations. Usually, a
  backward SPDE is solvable in the sense that there exists a
  diffusion term being considered as a part of the solution that
  helps to ensure that the solution is adapted to the driving Brownian
motions.
\par
There are also results for SPDEs with boundary conditions connecting the solution at different times, for instance, at initial time and at
terminal time. This category includes stationary type solutions for
forward SPDEs (see, e.g., Caraballo {\em et al } (2004),
 Chojnowska-Michalik (19987), Chojnowska-Michalik and Goldys
 (1995), Duan {\em et al} (2003), Mattingly (1999)
Mohammed  {\em et al} (2008), Sinai (1996), and the references here).  There are also results for periodic solutions of SPDEs
 (Chojnowska-Michalik (1990), Feng and Zhao (2012), Kl\"unger  (2001)).
 As was mentioned in Feng and Zhao (2012), it is difficult to expect that, in general, a SPDE has a periodic
 in time solution $u(\cdot,t)|_{t\in[0,T]}$ in a usual sense of exact equality  $u(\cdot,t)=u(\cdot,T)$ that holds almost surely
 given that $u(\cdot,t)$ is adapted to some Brownian motion.
The periodicity of the solutions of stochastic equations was
usually considered
 in the sense of the distributions. In Feng and Zhao (2012),
 the periodicity was established  in a  stronger sense as a "random
 periodic solution (see Definition 1.1 from Feng and Zhao (2012)).
\iindex{Dokuchaev (2012) considered  backward SPDEs with quite
general non-local and time and space boundary conditions. These
conditions cover  a setting where periodicity condition hold almost
surely, as well as more general conditions  $\kappa
u(\cdot,0)=u(\cdot,T)+\xi$ a.e.,where $\kappa\in[-1,1]$ and $\xi$ is
a random variable. Note that $u(\cdot,0)$ was assumed to be
non-random.} This was a novel setting comparing with the periodic
conditions for the distributions, or with conditions from Kl\"unger
(2001) and Feng and Zhao (2012), or with conditions for expectations
from Dokuchaev (2008).

The present paper addresses these and related problems again. We
consider linear  Dirichlet  condition at the boundary of the state
domain; the equations are of  a parabolic type  and are not
necessary self-adjoint. The standard boundary value Cauchy condition
at the one fixed time is replaces by a condition that mixes in one
equation the terminal value of the solution and a functional of the
 entire solution.  This setting covers conditions such as
 $\t^{-1}\int_0^\t u(\cdot,t)dt=u(\cdot,T)$ a.s., as well as  more general
conditions.

 We present sufficient conditions for existence and
regularity of the solutions  in $L_2$-setting (Theorem \ref{Th6}).
These results open a way to extend applications of backward SPDEs on
the problems with non-local in time space boundary conditions. Our approach is  based
 on the contraction mapping theorem in a $L_\infty$-space. 

A less  general case was
considered in Dokuchaev (2012b), where the boundary condition
was connecting $u(\cdot,T)$ with the expectations of the past values
of $u$.  In Dokuchaev (2012c), related forward and backward SPDEs were studied in an unified framework.  
In Dokuchaev (2012b,c), the approach was based
 on the Fredholm Theorem in a $L_2$-space;
 this approach is not applicable for  the setting  considered in the present paper.

\section{The problem setting and definitions}
 We
are given a standard  complete probability space $(\O,\F,\P)$ and a
right-continuous filtration $\F_t$ of complete $\s$-algebras of
events, $t\ge 0$. We assume that $\F_0$ is the $\P$-augmentation of
the set $\{\emptyset,\O\}$. We are given also a $N$-dimensional
Wiener process $w(t)$ with independent components;  it is a Wiener
process with respect to $\F_t$.
\par
Assume that we are given a bounded open domain $D\subset\R^n$  with
$C^2$-smooth boundary $\p D$. Let $T>0$ be given, and let $Q\defi
D\times [0,T]$. \par
 We will study the following boundary value
problem in $Q$
\begin{eqnarray} %4.1
\label{parab1} &&d_tu+(\A u+ \varphi)\,dt +\sum_{i=1}^N
B_i\chi_idt=\sum_{i=1}^N\chi_i(t)dw_i(t), \quad t\ge 0,
\\\label{parab10}
&& u(x,t,\o)\,|_{x\in \p D}=0
\\ &&u(\cdot, T)-\G u(\cdot)=\xi.
\label{parab2}%3.1
\end{eqnarray}
Here $u=u(x,t,\o)$, $\varphi=\varphi(x,t,\o)$, $\xi=\xi(x,\o)$,
$\chi_i=\chi_i(x,t,\o)$,
 $(x,t)\in Q$,   $\o\in\O$.\par
   In (\ref{parab2}), $\G$ is a linear operator that maps functions
defined on $Q\times \O$  to functions defines on $D\times \O$. For
instance, the case where  $\G u=u(\cdot,0)$ is not excluded; this
case corresponds to the periodic type boundary condition \baa
u(\cdot,T)-u(\cdot,0)=\xi.\label{period}\eaa
\par
  In (\ref{parab1}),  \baa \A
v=\sum_{i,j=1}^nb_{ij}(x,t,\o)\frac{\p^2 v}{\p x_i \p x_j}(x)
+\sum_{i=1}^n f_i(x,t,\o)\frac{\p v}{\p x_i
}(x)+\lambda(x,t,\o)v(x), \label{A}\eaa
  and \be\label{B}
B_iv\defi\frac{dv}{dx}\,(x)\,\beta_i(x,t,\o),\quad i=1,\ldots ,N.
\ee
\par
We assume that the functions $b(x,t,\o):
\R^n\times[0,T]\times\O\to\R^{n\times n}$, $\b_j(x,t,\o):
\R^n\times[0,T]\times\O\to\R^n$, $f(x,t,\o):
\R^n\times[0,T]\times\O\to\R^n$, $\lambda(x,t,\o):
\R^n\times[0,T]\times\O\to\R$,   $\chi_i(x,t,\o): \R^n\times
[0,T]\times\O\to\R$, and $\varphi (x,t,\o): \R^n\times
[0,T]\times\O\to\R$ are progressively measurable with respect to
$\F_t$ for all $x\in\R^n$, and the function $\xi(x,\o):
\R^n\times\O\to\R$ is $\F_0$-measurable for all $x\in\R^n$.

In fact,
we will also consider $\varphi$  from wider classes. In
particular, we will consider generalized functions $\varphi$.

We assume  $\lambda(x,t,\o)\le 0$ a.e., and
  $b_{ij}, f_i, x_i$ are the
components of $b$, $f$, and $x$ respectively.
\subsection*{Spaces and classes of functions} %2
We denote by $\|\cdot\|_{ X}$ the norm in a linear normed space
$X$, and
 $(\cdot, \cdot )_{ X}$ denote  the scalar product in  a Hilbert space $
X$.
\par
We introduce some spaces of real valued functions.
\par
 Let $G\subset \R^k$ be an open
domain, then ${W_q^m}(G)$ denote  the Sobolev  space of functions
that belong to $L_q(G)$ together with the distributional
derivatives up to the $m$th order, $q\ge 1$.
\par
 We denote  by $|\cdot|$ the Euclidean norm in $\R^k$, and $\bar G$ denote
the closure of a region $G\subset\R^k$.
\par Let $H^0\defi L_2(D)$,
and let $H^1\defi \stackrel{\scriptscriptstyle 0}{W_2^1}(D)$ be the
closure in the ${W}_2^1(D)$-norm of the set of all smooth functions
$u:D\to\R$ such that  $u|_{\p D}\equiv 0$. Let $H^2=W^2_2(D)\cap
H^1$ be the space equipped with the norm of $W_2^2(D)$. The spaces
$H^k$ and $W_2^k(D)$ are called  Sobolev spaces, they are Hilbert
spaces, and $H^k$ is a closed subspace of $W_2^k(D)$, $k=1,2$.
\par
 Let $H^{-1}$ be the dual space to $H^{1}$, with the
norm $\| \,\cdot\,\| _{H^{-1}}$ such that if $u \in H^{0}$ then
$\| u\|_{ H^{-1}}$ is the supremum of $(u,v)_{H^0}$ over all $v
\in H^1$ such that $\| v\|_{H^1} \le 1 $. $H^{-1}$ is a Hilbert
space.
\par We shall write $(u,v)_{H^0}$ for $u\in H^{-1}$
and $v\in H^1$, meaning the obvious extension of the bilinear form
from $u\in H^{0}$ and $v\in H^1$.
\par
We denote by $\oo\ell _{k}$ the Lebesgue measure in $\R^k$, and we
denote by $ \oo{{\cal B}}_{k}$ the $\sigma$-algebra of Lebesgue
sets in $\R^k$.
\par
We denote by $\oo{{\cal P}}$  the completion (with respect to the
measure $\oo\ell_1\times\P$) of the $\s$-algebra of subsets of
$[0,T]\times\O$, generated by functions that are progressively
measurable with respect to $\F_t$.
\par
 We  introduce the spaces
 \baaa
 &&X^{k}(s,t)\defi L^{2}\bigl([ s,t ]\times\Omega,
{\oo{\cal P }},\oo\ell_{1}\times\P;  H^{k}\bigr), \quad\\ &&Z^k_t
\defi L^2\bigl(\Omega,{\cal F}_t,\P; H^k\bigr),\\
&&\CC^{k}(s,t)\defi C\left([s,t]; Z^k_T\right), \qquad k=-1,0,1,2,
\\&& \X^k_c= L^{2}\bigl([ 0,T ]\times\O,\, \oo{{\cal P}
},\oo\ell_{1}\times\P;\; C^k(\oo D)\bigr),\quad k\ge 0. \eaaa
%$\V^{k}(s,T)=L^{2}\bigl([s,T ],\oo\ell_{1}\times\P;W^{k}_2(D)\bigr)$,
%$ k=0,2,..$,
  The
spaces $X^k(s,t)$ and $Z_t^k$  are Hilbert spaces.
 \par
We introduce the spaces $$ Y^{k}(s,t)\defi
X^{k}(s,t)\!\cap \CC^{k-1}(s,t), \quad k=1,2, $$ with the norm $ \|
u\| _{Y^k(s,T)}
\defi \| u\| _{{X}^k(s,t)} +\| u\| _{\CC^{k-1}(s,t)}. $
For brevity, we shall use the notations
 $X^k\defi X^k(0,T)$, $\CC^k\defi \CC^k(0,T)$,
and  $Y^k\defi Y^k(0,T)$.

We also introduce spaces $\CC^k_{PC}$ consisting of $u\in \CC^k$
such that either $u\in \CC^k$ or there exists $\t=\t(u)\in [0,T]$
such that $\|u(\cdot,t)\|_{Z_T^k}$ is bounded,  $u(\cdot,t)$ is
continuous in $Z_T^k$ in $t\in[0,\t]$, and
 $u(\cdot,t)$ is continuous in $Z_T^k$ in $t\in[\t+\e,T]$ for any $\e>0$.

Finally, we introduce the spaces
 \baaa
 &&\W\defi L^{\infty}\bigl([ 0,T]\times\Omega,
{\oo{\cal P }},\oo\ell_{1}\times\P;  L_{\infty}(D)\bigr)\cap \CC_{PC}^0(0,T), \quad\\ &&\V
\defi L^{\infty}\bigl(\Omega,{\cal F}_T,\P;  L_{\infty}(D) \bigr).
\eaaa
\index{
We denote $\W=\W(0,T)$.}\par

\subsection*{Conditions for the coefficients}
 To proceed further, we assume that Conditions
\ref{cond3.1.A}-\ref{condK} remain in force throughout this paper.
 \begin{condition} \label{cond3.1.A} The matrix  $b=b^\top$ is
symmetric  and bounded. In addition, there exists a constant
$\d>0$ such that
\be
 \label{Main1} y^\top  b
(x,t,\o)\,y-\frac{1}{2}\sum_{i=1}^N |y^\top\b_i(x,t,\o)|^2 \ge
\d|y|^2 \quad\forall\, y\in \R^n,\ (x,t)\in  D\times [0,T],\
\o\in\O. \ee
\end{condition}
\begin{condition}\label{cond3.1.B}
The functions  $f(x,t,\o)$, $\lambda (x,t,\o)$, and $\b_i(x,t,\o)$ and
are bounded. These functions  are differentiable
in $x$ for a.e. $t,\o$, and the corresponding derivatives are
bounded. In addition, $b\in \X_c^3$, $\w f\in\X_c^2$,
$\lambda\in\X^1_c$, $ \b_i\in\X_c^3$, and $\b_i(x,t,\o)=0$ for
$x\in \p D$, $i=1,...,N$.
\end{condition}
Let $\Ind$ denote  the indicator function
\begin{condition}\label{condK} The mapping $\G:
\W\to \V$ is linear and continuous and such that $\|\G
u\|_{\V}\le \|u\|_{\W}$ for any $u\in\W$, and that there exists
$\t<T$ such that $\G u=\G (\Ind_{\{t\le \t\}}u)$.
\end{condition}
\begin{example} {\rm
Condition \ref{condK} is satisfied for the following operators:
\begin{enumerate}
\item $\G u=\kappa u(\cdot,0)$, $\kappa\in[-1,1]$;
\item
$ (\G u)(x,\o)=\kappa u(x,t_1,\o),\quad t_1\in[0,T);$
\item
$(\G u)(x,\o)=\zeta(\o) u(x,t_1,\o),\quad t_1\in[0,T),\qquad
\zeta\in L_{\infty}(\O,\P,\F_T,\P),\quad |\zeta(\o)|\le 1\quad
\hbox{a.s.} $;
\item
$ (\G u)(x,\o)=\a_1 u(x,t_1,\o)+\a_2 u(x,t_2,\o),\quad
t_1,t_2\in[0,T),\quad |\a_1|+|\a_2|\le 1$;
\item \baaa (\G
u)(x,\o)=\int_0^{\t}k(t)u(x,t,\o)dt,\quad \t\in[0,T),\qquad
k(\cdot)\in L_{\infty}(0,\t),\qquad  \int_0^\t |k(t)|dt\le 1; \eaaa
  \item
\baaa (\G u)(x,\o)=\int_0^{\t}dt\int_Dk(t,y,x,\o)u(y,t,\o)dy, \eaaa
where $\t\in[0,T)$, $k(\cdot):[0,\t]\times D\times D\times\O$ is a
bounded measurable function from $L^{\infty}
(\O,\F_T,\P,L_{\infty}([0,\t]\times D\times D))$
 such that \baaa \esssup_{(x,\o)\in D\times\O}\int_0^\t dt\int_D |k(t,x,y,\o)|dy\le 1. \eaaa
\end{enumerate}
Convex combinations of operators from this list are also covered.}
\end{example}

Sometimes we shall omit $\o$.
\subsection*{The definition of solution}
\begin{proposition} %2.2
\label{propL} Let $\zeta\in X^0$,
 let a sequence  $\{\zeta_k\}_{k=1}^{+\infty}\subset
L^{\infty}([0,T]\times\O, \ell_1\times\P;\,C(D))$ be such that all
$\zeta_k(\cdot,t,\o)$ are progressively measurable with respect to
$\F_t$, and let $\|\zeta-\zeta_k\|_{X^0}\to 0$. Let $t\in [0,T]$ and
$j\in\{1,\ldots, N\}$ be given.
 Then the sequence of the
integrals $\int_0^t\zeta_k(x,s,\o)\,dw_j(s)$ converges in $Z_t^0$ as
$k\to\infty$, and its limit depends on $\zeta$, but does not depend
on $\{\zeta_k\}$.
\end{proposition}
\par
{\it Proof} follows from completeness of  $X^0$ and from the
equality
\begin{eqnarray*}
\E\int_0^t\|\zeta_{k}(\cdot,s,\o)-\zeta_m(\cdot,s,\o)\|_{H^0}^2\,ds
=\int_D\,dx\,\E\left(\int_0^t\big(\zeta_k(x,s,\o)-
\zeta_m(x,s,\o)\big)\,dw_j(s)\right)^2.
\end{eqnarray*}
\begin{definition} %{2.1}
\rm Let $\zeta\in X^0$, $t\in [0,T]$, $j\in\{1,\ldots, N\}$, then we
define $\int_0^t\zeta(x,s,\o)\,dw_j(s)$ as the limit  in $Z_t^0$ as
$k\to\infty$ of a sequence $\int_0^t\zeta_k(x,s,\o)\,dw_j(s)$, where
the sequence $\{\zeta_k\}$ is such  as in Proposition \ref{propL}.
\end{definition}
\begin{definition} %3.1
\label{defsolltion} \rm Let $u\in Y^1$, $\chi_i\in X^0$,
$i=1,...,N$, and $\varphi\in X^{-1}$. We say that equations
(\ref{parab1})-(\ref{parab10}) are satisfied if \baaa
&&u(\cdot,t,\o)=u(\cdot,T,\o)+ \int_t^T\big(\A u(\cdot,s,\o)+
\varphi(\cdot,s,\o)\big)\,ds \ \nonumber
\\&&\hphantom{xxx}+ \sum_{i=1}^N
\int_t^TB_i\chi_i(\cdot,s,\o)ds-\sum_{i=1}^N
\int_t^T\chi_i(\cdot,s)\,dw_i(s)
%\eqno(3.2)
\label{intur} \eaaa for all $r,t$ such that $0\le r<t\le T$, and
this equality is satisfied as an equality in $Z_T^{-1}$.
\end{definition}
Note that the condition on $\p D$ is satisfied in the  sense that
$u(\cdot,t,\o)\in H^1$ for a.e. \ $t,\o$. Further, $u\in Y^1$, and
the value of  $u(\cdot,t,\o)$ is uniquely defined in $Z_T^0$ given
$t$, by the definitions of the corresponding spaces. The integrals
with $dw_i$ in (\ref{intur}) are defined as elements of $Z_T^0$. The
integral with $ds$ in (\ref{intur}) is defined as an element of
$Z_T^{-1}$. In fact, Definition \ref{defsolltion} requires for
(\ref{parab1}) that this integral must be equal  to an element of
$Z_T^{0}$ in the sense of equality in $Z_T^{-1}$.
\section{The main results}%{Problems with non-local conditions}
\begin{theorem}
\label{Th6}  Problem (\ref{parab1})-(\ref{parab2}) has
 a unique solution
$(u,\chi_1,...,\chi_N)$ in the class $Y^1\times(X^0)^N$ for any
$\varphi\in \W$ and $\xi\in  Z^0_T$. This solution is such that
$u\in \W$. In addition, \baa \label{3.Th6} \|u\|_{\W}+\| u
\|_{Y^1}+\sum_{i=1}^N\|\chi_i\|_{X^0}\le C \left(\| \varphi \|
_{\W}+\|\xi\|_{\V}\right), \eaa  where $C>0$
 does not depend on $\varphi$ and  $\xi$.
\end{theorem}
\index{ and let  $\w\lambda(x,t,\o)\le c_\lambda$ a.e., and where
$c_\lambda<0$ is given.}
\section{Proofs}
%\subsection*{Problems for forward equations}
Let $s\in (0,T]$, $\varphi\in X^{-1}$ and $\Phi\in Z^0_s$. Consider
the problem \be \label{4.1}
\begin{array}{ll}
d_tu+\left( \A u+ \varphi\right)dt +
\sum_{i=1}^NB_i\chi_i(t)dt=\sum_{i=1}^N\chi_i(t)dw_i(t), \quad t\le s,\\
u(x,t,\o)|_{x\in \p D}, \\
 u(x,s,\o)=\Phi(x,\o).
\end{array}
 \ee
 \par
The following lemma represents an analog of the so-called "the first
energy inequality", or "the first fundamental inequality" known for
deterministic parabolic equations (see, e.g., inequality (3.14) from
Ladyzhenskaya (1985), Chapter III).
\begin{lemma}
\label{lemma1} Assume that Conditions \ref{cond3.1.A}--\ref{condK}
are satisfied.  Then problem (\ref{4.1}) has an unique solution a
unique solution $(u,\chi_1,...,\chi_N)$ in the class
$Y^1\times(X^0)^N$  for any $\varphi\in X^{-1}(0,s)$, $\Phi\in
Z_s^0$, and \be \label{4.2} \| u
\|_{Y^1(0,s)}+\sum_{i=1}^N\|\chi_i\|_{X^0}\le C \left(\| \varphi \|
_{X^{-1}(0,s)}+\|\Phi\|_{Z^0_s}\right), \ee where $C>0$
does not depend on $\varphi$ and $\xi$.
\end{lemma}
(See, e.g., Dokuchaev (1991) or Theorem 4.2 from Dokuchaev (2010)).
\par Note that the solution $u=u(\cdot,t)$
is continuous in $t$ in $L_2(\O,\F,\P,H^0)$, since
$Y^1(0,s)=X^{1}(0,s)\!\cap \CC^{0}(0,s)$.
\par
Introduce  operators $L_s:X^{-1}(0,s)\to Y^1(0,s)$ and
$\L_s:Z^0_s\to Y^1(0,s)$, such that $u=L_s\varphi+\L_s\Phi,$ where
 $(u,\chi_1,...,\chi_N)$ is the solution of  problem (\ref{4.1})  in the class
$Y^2\times(X^1)^N$. By Lemma \ref{lemma1}, these linear operators
are continuous.
\par
Introduce   operators $\Q:Z_T^0\to Z_T^0$ and $\TT:X^{-1}\to Z_T^0$
such that $\Q \Phi=\G\L_T\Phi$ and $\TT\varphi=\G L_T\varphi$, i.e,  $\Q\Phi+\TT\varphi= \G u$,
    where
$u$ is the solution in $Y^1$ of   problem (\ref{4.1}) with $s=T$,
$\varphi\in X^{-1}$, and $\Phi\in Z_T^0$.

It is easy to see that if the operator $\G:Y^1\to Z_T^0$ is continuous, then
the operators   $\Q:Z_T^0\to Z_T^0$ and $\TT:X^{-1}\to Z_T^0$ are linear and continuous. In particular,
$\|\Q\|\le \|\G\|\|\L_T\|$, where $\|\Q\|$, $\|\G\|$, and
$\|\L_T\|$, are the norms of the operators $\Q: Z_T^0\to Z_T^0$,
$\G: Y^1\to Z_T^0$, and $\L_T: Z_T^0\to Y^1$, respectively.

\begin{lemma}\label{lemmaQ} Assume that the operator $\G:Y^1\to Z_T^0$ is continuous. If the operator $(I-\Q)^{-1}:Z_T^0\to Z_T^0$ is also continuous
then problem (\ref{4.1})   has a unique solution
$(u,\chi_1,...,\chi_N)$ in the class $Y^1\times(X^0)^N$ for any
$\varphi \in X^{-1}$, $\Phi \in Z_T^0$. For this solution, \baa
\label{Q} u=L_T\varphi+\L_T (I-\Q)^{-1}(\xi+\TT\varphi) \eaa
 and \baaa
\label{Qes} \| u \|_{Y^1(0,s)}+\sum_{i=1}^N\|\chi_i\|_{X^0}\le C
\left(\| \varphi \| _{X^{-1}(0,s)}+\|\Phi\|_{Z^0_s}\right), \eaaa
where $C=C({\cal P})$ does not depend on $\varphi$ and $\xi$.
\end{lemma}
\par
{\it Proof of Lemma \ref{lemmaQ}}.    For brevity, we denote
$u(\cdot,t)=u(x,t,\o)$. Clearly,
 $u\in Y^1$ is the solution of   problem
(\ref{parab1})-(\ref{parab2}) with some $(\chi_1,...,\chi_N)\in
(X^0)^N$ if and only if \baa
&&u=\L_Tu(\cdot,T)+L_T\varphi, \label{Q1}\\
&&u(\cdot,T)-\G u=\xi. \label{Q2}\eaa
 Since $ \G u=\Q
u(\cdot,T)+\TT\varphi$, equation (\ref{Q2})  can be rewritten as
\baa
 u(\cdot,T)-\Q u(\cdot,T)-\TT\varphi
=\xi. \label{Q3}\eaa  By the continuity of $(I-\Q)^{-1}$, equation
(\ref{Q3}) can be rewritten as
$$ u(\cdot,T)=(I-\Q)^{-1}(\xi +\TT\varphi). $$  Therefore, equations (\ref{Q1})-(\ref{Q2}) imply that \baaa \label{4.3}
u=L_T\varphi+\L_T u(\cdot,T)=L_T\varphi+\L_T
(I-\Q)^{-1}(\xi+\TT\varphi). \eaaa Further, let us show that if
(\ref{Q}) holds then equations (\ref{Q1})-(\ref{Q2})  hold. Let $u$
be defined by (\ref{Q}). Since $u=L_T\varphi+\L_T u(\cdot,T)$, it
follows that $u(\cdot,T)=(I-\Q)^{-1}(\xi+\TT\varphi)$. Hence \baaa
u(\cdot,T)-\Q u(\cdot,T)=\xi+\TT\varphi,\eaaa i.e., $u(\cdot,T)-\G
\L_T u(\cdot,T)=\xi+\TT\varphi=\xi+\G L_T\varphi.$   Hence \baaa
u(\cdot,T)-\G [\L_T u(\cdot,T)+L_T\varphi]=\xi.\eaaa This means that
(\ref{Q1})-(\ref{Q2}) hold. Then the proof of Lemma \ref{lemmaQ}
follows. $\Box$
\par
Let functions ${\ww\b_i: Q\times \O \to \R^n}$, $i=1,\ldots, M$, be
such that $$ 2b(x,t,\o)=\sum_{i=1}^N\b_i(x,t,\o)\,\b_i(x,t,\o)^\top
+\sum_{j=1}^M\,\ww\b_j(x,t,\o)\,\ww\b_j(x,t,\o)^\top, $$ and $\ww
\b_i$ has the similar properties as $\b_i$. (Note that, by Condition
\ref{cond3.1.A}, $2b>\sum_{i=1}^N\b_i\b_i^\top$).
\par
 Let
$\ww w(t)=(\ww w_1(t),\ldots, \ww w_M(t))$ be a new Wiener process
independent on $w(t)$. Let $a\in L_2(\O,\F,\P;\R^n)$ be a vector
such that $a\in D$. We assume also that $a$ is independent from
$(w(t)-w(t_1),\w w(t)-\w w(t_1))$ for all $t>t_1>s$. Let $s\in[0,T)$
be given. Consider the following Ito equation
\begin{eqnarray}
%\begin{array}{c}
\label{yxs} &&dy(t) =
f(y(t),t)\,dt+\sum_{i=1}^N\b_i(y(t),t)\,dw_i(t) +\sum_{j=1}^M\ww
\b_j(y(t),t)\,d \ww w_j(t),
%\vspace*{20pt}
\nonumber\\ [-6pt] &&y(s)=x.
%\eqno{3.12}
%\end{array}
\end{eqnarray}
\par
Let  $y(t)=y^{a,s}(t)$ be the solution of (\ref{yxs}), and let
$\tau^{a,s}\defi\inf\{t\ge s:\ y^{a,s}(t)\notin D\}$.
\begin{lemma}\label{propnu} For any $\vartheta>0$, there exists  $\nu=\nu(\vartheta)\in(0,1)$
that depends only on $D,\A,B_j$ and such that $\P_s(\tau^{x,s}>
s+\vartheta)\le \nu$ a.s. for all $s\ge 0$, and for any $x\in D$.
\end{lemma}
\par
Note that if the functions $ f(x,t,\o)= f(x)$ and $\b(x,t,\o)=\b(x)$
are non-random and constant in $t$, then existence of $\nu\in(0,1)$
such that $\P(\tau^{a,s}> s+\vartheta )\le \nu$ $(\forall a,s)$ is obvious.
\par {\it Proof of Lemma
\ref{propnu}.}  In this proof, we will follow the approach from
Dokuchaev (2004), p.296. Let $\mu =(\w f,\b,x,s)$.\par Clearly, there
exists a finite interval $D_1\defi (d_1,d_2)\subset \R$ and a
bounded domain $D_{n-1}\subset \R^{n-1}$ such that $D\subset
D_1\times D_{n-1}$.
\par
For $(x,s)\in D\times [0,T)$, let $\tau_1^{x,s}\defi\inf\{t\ge s: \
y^{x,s}_1(t)\notin D_1\}$, where $y^{x,s}_1(t)$ is the first
component of the vector
$y^{x,s}(t)=(y^{x,s}_1(t),...,y^{x,s}_n(t))$. We have that
\be\label{dd1} \P_s(\tau^{x,s}>s+\vartheta )\le
\P_s(\tau_1^{x,s}>s+\vartheta )=\P_s(y_1^{x,s}(t)\in D_1\ \forall
t\in[s,s+\vartheta ]). \ee
\par
Let \baaa M^{\mu}(t)\defi
\sum_{k=1}^N\int_s^th_k(y^{x,s}(r),r)dw_i(r)+\sum_{k=N+1}^{N+M}\int_s^th_k(y^{x,s}(r),r)d\ww
w_i(r),\quad t\ge s, \eaaa where $h=(h_1,..,h_{N+M})$ is a vector
  that represents the first
row of the matrix \baaa (\b_1,...,\b_N,\w\b_1,...,\w\b_M) \eaaa with
the values in $\R^{n\times (N+M)}$.
\par
 Let $\w D_1\defi (d_1+K_1,d_2+K_2)$, where
$K_1\defi -d_2-\vartheta \sup_{x,t,\o}|\w f_1(x,t,\o)|$, $K_2\defi
-d_1+\vartheta \sup_{x,t}|\w f_1(x,t,\o)|$. Clearly, $\w D_1$
depends only on $n,D$, and $c_f$. It is easy to see that
\be\label{dd2} \P_s(y_1^{x,s}(t)\in D_1\ \forall t\in[s,s+\vartheta
])\le \P_s(M^{\mu}(t)\in \w D_1\ \forall t\in[s,s+\vartheta ]). \ee
Further, \be\label{dd3}h(y^{x,s}(t),t)^\top h(y^{x,s}(t),t)=
|h(y^{x,s}(t),t)|^2\in [\d,c_{\b}],\ee where \baaa \d=\inf_{x,s,\o,\
\xi\in\R^n:\ |\xi|=1}2\xi^\top b(x,t,\o)\xi,\quad
c_\b=\sup_{x,s,\o,\ \xi\in\R^n:\ |\xi|=1}2\xi^\top
b(x,t,\o)\xi.\eaaa Clearly, $M^{\mu}(t)$ is a martingale vanishing
at $s$ conditionally given $\F_s$ with quadratic variation process
$$[M^{\mu}]_t\defi \int_s^t|h(y^{x,s}(r),r)|^2dr,\qquad t\ge s.
$$
\par Let $\t^{\mu}(t)\defi \inf\{r\ge s:\ [M^{\mu}]_r>t-s\}$.
 Note that $\t^{\mu}(s)=s$, and
the function $\t^{\mu}(t)$ is  strictly increasing in $t>s$ given
$(x,s)$.  By Dambis--Dubins--Schwarz Theorem (see, e.g., Revuz and
Yor (1999)), the process $B^{\mu}(t)\defi M(\t^{\mu}(t))$ is a
Brownian motion  conditionally given $\F_s$ vanishing at $s$, i.e.,
$B^{\mu}(s)=0$, and $M^{\mu}(t)=B^{\mu}(s+[M^{\mu}]_t)$. Clearly,
\be\label{dd4} \ba \P_s(M^{\mu}(t)\in \w D_1\ \hphantom{x}\forall
t\in[s,s+\vartheta ])&=\P_s(B^{\mu}(s+[M^{\mu}]_t)\in \w D_1\
\hphantom{x}\forall
t\in[s,s+\vartheta ])\\
&\le \P_s(B^{\mu}(r)\in \w D_1\ \hphantom{x}\forall
r\in[s,s+[M^{\mu}]_{s+\vartheta }]). \ea \ee By (\ref{dd3}),
$[M^{\mu}]_{s+\vartheta }\ge \d \vartheta $ a.s. for all $x,s$.
Hence \be\label{dd5}\P_s(B^{\mu}(r)\in \w D_1\ \hphantom{x}\forall
r\in[s,s+[M^{\mu}]_{s+\vartheta }])\le\P_s(B^{\mu}(r)\in \w D_1\
\hphantom{x}\forall r\in[s,s+\d  \vartheta ]). \ee By
(\ref{dd1})--(\ref{dd2}) and (\ref{dd4})--(\ref{dd5}), it follows
that
$$
\ba \sup_{\mu}\P_s(\tau^{x,s}>s+\vartheta )\le \nu\defi \sup_{\mu}
\P_s(B^{\mu}(r)\in \w D_1\ \hphantom{x}\forall r\in[s,s+\d \vartheta
]), \ea
$$
and $\nu=\nu({\cal P})\in(0,1)$.
 This completes the
proof of Lemma \ref{propnu}.
 $\Box$
\par
 {\em Proof of Theorem \ref{Th6}}. For $t\ge s$,
set
\begin{eqnarray*}
\g^{a,s}(t) \defi\exp\left(-\int_s^t
\lambda(y^{a,s}(t),t)\,dt\right).
\end{eqnarray*}
\par Let $\Phi\in \V$ and $\varphi\in
\W$ be bounded.  By Theorem 4.1 from Dokuchaev (2011) again, we have
that, for any $s\in[0,T)$ and $u=L_T\xi+\L_T\Phi$, $u(\cdot,s)$ can be represented as \baa
u(x,s,\o)=\E\left\{\g^{x,s}(T)\Phi(y^{x,s}(T))\Ind_{\{ \tau^{x,s}\ge
T\}} +\int_s^{\tau^{x,s}}
\g^{x,s}(t)\,\varphi(y^{x,s}(t),t,\o)\,dt\,\bigr|\,\F_s\right\}.\hphantom{X}\label{repR}\eaa
This equality holds in $Z_s^0$ and for a.e. $x,\o$.
 It follows that
  \baa
\sup_{s\in[0,T]}\|u(\cdot,s)\|_{\V}\le \|\Phi\|_{\V}+
T\|\varphi\|_{\W}.\label{estRs}\eaa Hence \baa  \|\L_T\Phi\|_{\W}\le
\|\Phi\|_{\V},\quad \|L_T\varphi\|_{\W}\le
T\|\varphi\|_{\W}.\label{estR}\eaa By the assumptions on $\G$, it
follows that $\|\G u\|_{\V}\le \|u\|_{\W}.$ It follows that the
operators $\Q=\G\L_T:\V\to \V$ and  $\TT:\W\to \V$ are  bounded. Let
$\|\Q\|_{\V,\V}$ be the norm of the operator $\Q:\V\to \V$.
\par
It follows from (\ref{estRs})  and from the properties of  $\G$ that
$\|\Q\|_{\V,\V}\le 1$. Let us refine this estimate.
\par \index{Assume first that $c_\lambda=0$.that $\w\lambda$ is replaced
by $\ww\lambda=\w\lambda+|c_\lambda|$.}

\par
Let $u=\L_T\Phi$, $s\in[0,T]$. Let $y(t)=y^{x,s}(t)$ be
the solution of  Ito equation (\ref{yxs}) with the initial condition
$y(s)=x$.
 For the brevity, we will use notations $\P_s(\cdot)\defi
\P(\cdot|\F_s)$ and  $\E_s(\cdot)\defi \E(\cdot|\F_s)$.
 By
(\ref{repR}), it follows that \baaa
&&\|u(\cdot,s)\|_{\V}=\esssup_{x,\o}\E_s\g^{x,s}(T)\Phi(y^{x,s}(T))\Ind_{\{
\tau^{x,s}\ge T\}}\\  &&\le \esssup_{x,\o} \left[\E_s\Ind_{\{
\tau^{x,s}\ge T\}}^2\right]^{1/2} \esssup_{x,\o}
\left[\E_s\Phi(y^{x,s}(T))^2\right]^{1/2}
\\
 &&\le   \esssup_{x,\o} \left[\E_s\Ind_{\{ \tau^{x,s}\ge
T\}}^2\right]^{1/2} \|\Phi\|_{\V}=\esssup_{x,\o} \P_s(\tau^{x,s}\ge
T) ^{1/2} \|\Phi\|_{\V}. \eaaa If $s<\t$ then $\{\tau^{x,s}\ge
T\}\subseteq \{\tau^{x,s}\ge s+\vartheta\}$, where $\vartheta\defi
T-\t>0$. Hence
 \baaa \|u(\cdot,s)\|_{\V}\le\esssup_{x,\o}
\P_s(\tau^{x,s}\ge s+\vartheta) ^{1/2} \|\Phi\|_{\V},\quad s\le \t.
\eaaa
 By Lemma
\ref{propnu}, it follows that there exists $\nu=\nu(\vartheta,{\cal
P})\in (0,1)$ such that $\P_s(\tau^{x,s}\ge s+\vartheta)<\nu$ a.s.
It follows that \baaa \|u(\cdot,s)\|_{\V} \le \nu^{1/2}
\|\Phi\|_{\V},\quad s\le \t \eaaa and  \baaa \|\Ind_{\{s\le
\t\}}u\|_{\W} \le \nu^{1/2} \|\Phi\|_{\V}. \eaaa By the assumptions
on $\G$, it follows that \baaa \|\G u\|_{\V}=\|\G (\Ind_{\{s\le
\t\}}u)\|_{\V} \le \nu^{1/2} \|\Phi\|_{\V},\quad s\le \t. \eaaa
 It follows that $\|\Q\|_{\V,\V}\le  \nu^{1/2}<1$.
Hence  the operator
    $(I-\Q)^{-1}:\V\to \V$ is bounded. Let \baa
\label{uQ} u=L_T\varphi+\L_T (I-\Q)^{-1}(\xi+\TT\varphi).\eaa By the
assumptions on $\G$ and by (\ref{repR})-(\ref{estR}), it follows
that  $\xi+\TT\varphi=\xi+\G L_T\varphi\in\V\subset Z_T^0$. Hence
$(I-\Q)^{-1}(\xi+\TT\varphi)\in\V\subset Z_T^0$. By the properties
of $\L_T$ and $L_T$, it follows that $u\in Y^1$. By
(\ref{repR})-(\ref{estR}) again, it follows that $u\in\W$.
 Similarly to the proof of Lemma \ref{lemmaQ}, it can be shown that $u$ is a part
 of the unique  solution
$(u,\chi_1,...,\chi_N)\in Y^1\times(X^0)^N$ of problem
(\ref{parab1})-(\ref{parab2}). Estimate (\ref{3.Th6}) follows from
the continuity of the corresponding operators in (\ref{uQ}). Then
the proof of Theorem \ref{Th6} follows. $\Box$
\subsection*{Acknowledgment} This work  was
supported by ARC grant of Australia DP120100928 to the
author.
\section*{References} $\hphantom{XX}$
Al\'os, E., Le\'on, J.A., Nualart, D. (1999).
 Stochastic heat equation with random coefficients
 {\it
Probability Theory and Related Fields} {\bf 115} (1), 41--94.
\par
Bally, V., Gyongy, I., Pardoux, E. (1994). White noise driven
parabolic SPDEs with measurable drift. {\it Journal of Functional
Analysis} {\bf 120}, 484--510.
\par Caraballo,T., P.E. Kloeden, P.E.,
Schmalfuss, B. (2004). Exponentially stable stationary solutions for
stochastic evolution equations and their perturbation, Appl. Math.
Optim. {\bf  50}, 183--207.
\par
 Chojnowska-Michalik, A. (1987). On processes of Ornstein-Uhlenbeck type in
Hilbert space, Stochastics 21, 251--286.

\par
Chojnowska-Michalik, A. (1990). Periodic distributions for linear
equations with general additive noise, Bull. Pol. Acad. Sci. Math.
38 (1–12) 23--33.

\par
Chojnowska-Michalik, A., and Goldys, B. (1995). {Existence,
uniqueness and invariant measures for stochastic semilinear
equations in Hilbert spaces},  {\it Probability Theory and Related
Fields},  {\bf 102}, No. 3, 331--356.
\par
Da Prato, G., and Tubaro, L. (1996). { Fully nonlinear stochastic
partial differential equations}, {\it SIAM Journal on Mathematical
Analysis} {\bf 27}, No. 1, 40--55.
\par
Dokuchaev, N.G. (1992). { Boundary value problems for functionals
of
 Ito processes,} {\it Theory of Probability and its Applications}
 {\bf 36} (3), 459-476.
 \par Dokuchaev, N.G. (2004).
Estimates for distances between first exit times via parabolic
equations in unbounded cylinders. {\it Probability Theory and
Related Fields}, {\bf 129}, 290 - 314.
\par
Dokuchaev, N.G. (2005).  Parabolic Ito equations and second
fundamental inequality.  {\it Stochastics} {\bf 77} (2005), iss. 4.,
pp. 349-370.
\par
 Dokuchaev N. (2008) Parabolic Ito equations with mixed in time
conditions.
{\it Stochastic Analysis and Applications} {\bf 26}, Iss. 3, 562--576. %May 2008
\par
Dokuchaev, N. (2010). Duality and semi-group property for backward
parabolic Ito equations. {\em Random Operators and Stochastic
Equations. } {\bf 18}, 51-72.
\par
Dokuchaev, N. (2011). Representation of functionals  of Ito
processes in bounded domains. {\em Stochastics} {\bf 83}, No. 1,
45--66.
\par
 Dokuchaev, N. (2012a).
Backward parabolic Ito equations and second fundamental inequality.
{\em Random Operators and Stochastic Equations} {\bf 20}, iss.1,
69-102.
\par
 Dokuchaev, N. (2012b). On almost surely
periodic and almost periodic solutions of backward SPDEs.
Working paper: http://arxiv.org/abs/1208.5538.
\par
Dokuchaev, N. (2012c). On forward and backward SPDEs with
non-local boundary conditions.  Working paper in  arXiv (submitted).  
\par
Du K., and Tang, S. (2012). Strong solution of backward stochastic
partial differential equations in $C^2$ domains.  {\em  Probability
Theory and Related Fields)} {\bf 154}, 255--285.
\par
Duan J.,  Lu K., Schmalfuss B. (2003). Invariant manifolds for
stochastic partial differential equations.{\em Ann. Probab.}{\bf 31}
 2109–2135.
\par
Feng C., Zhao H. (2012). Random periodic solutions of SPDEs via
integral equations and Wiener-Sobolev  compact embedding. {\em
Journal of Functional Analysis} {\bf 262}, 4377--4422.
\par
Gy\"ongy, I. (1998). Existence and uniqueness results for semilinear
stochastic partial differential equations. {\it Stochastic Processes
and their Applications} {\bf 73} (2), 271-299.
\par
Kl\"unger, M. (2001). Periodicity and Sharkovsky's theorem for
random dynamical systems, {\em Stochastic and Dynamics} {\bf 1},
iss.3, 299--338.
\par Krylov, N. V. (1999). An
analytic approach to SPDEs. Stochastic partial differential
equations: six perspectives, 185--242, Mathematical Surveys and
Monographs, {\bf 64}, AMS., Providence, RI, pp.185-242.
\par
Ladyzhenskaia, O.A. (1985). {\it The Boundary Value Problems of
Mathematical Physics}. New York: Springer-Verlag.
\par
Liu, Y., Zhao, H.Z (2009). Representation of pathwise stationary
solutions of stochastic Burgers equations, {\em Stochactics and
Dynamics} {\bf  9} (4), 613--634.
\par
Maslowski, B. (1995). { Stability of semilinear equations with
boundary and pointwise noise}, {\it Annali della Scuola Normale
Superiore di Pisa - Classe di Scienze} (4), {\bf 22}, No. 1,
55--93.
\par
Mattingly. J. (1999). Ergodicity of 2D Navier–Stokes equations with random forcing and large viscosity. {\em Comm. Math.
Phys.} 206 (2),  273–288.
\par
Mohammed S.-E.A., Zhang T.,  Zhao H.Z. (2008). The stable manifold theorem for semilinear stochastic evolution equations
and stochastic partial differential equations. {\em Mem. Amer. Math. Soc.} 196 (917)  1–105.
\par
Pardoux, E. (1993). Bulletin des Sciences Mathematiques, 2e Serie,
 {\bf 117}, 29-47.
\par
 Revuz, D., and Yor, M. (1999). {\it Continuous Martingales
and Brownian Motion}. Springer-Verlag: New York.
\par
Rodkina, A.E. (1992). On solutions of stochastic equations with
almost surely periodic
    trajectories.  {\em Differ. Uravn}. 28, No.3, 534--536 (in Russian).
\par
Rozovskii, B.L. (1990). {\it Stochastic Evolution Systems; Linear
Theory and Applications to Non-Linear Filtering.} Kluwer Academic
Publishers. Dordrecht-Boston-London.
%315 p.
\par
Sinai, Ya. (1996). Burgers system driven by a periodic stochastic
flows, in: Ito's Stochastic Calculus and Probability Theory,
Springer, Tokyo, 1996, pp. 347–353.
\par
Walsh, J.B. (1986). An introduction to stochastic partial
differential equations, Lecture Notes in Mathematics {\bf 1180},
Springer Verlag.
\par
Yong, J., and Zhou, X.Y. (1999). { Stochastic controls: Hamiltonian
systems and HJB equations}. New York: Springer-Verlag.
\par
 Zhou, X.Y. (1992). { A duality analysis on stochastic partial
differential equations}, {\it Journal of Functional Analysis} {\bf
103}, No. 2, 275--293.
\end{document}